\documentclass[12pt]{article}
\usepackage{amsmath,amsthm,amssymb,amsfonts}

\newfont{\lie}{eufm10 at 11pt}
\newfont{\liepequenos}{eufm10 at 10pt}

\newfont{\corpos}{msbm10 at 11pt}
\newfont{\corpospequenos}{msbm10 at 10pt}

\newcommand{\C}{\mbox{\corpos \symbol{67}}}          
\newcommand{\Hamil}{\mbox{\corpos \symbol{72}}}      
\newcommand{\Octoni}{\mbox{\corpos \symbol{79}}}     
\newcommand{\Pro}{\mbox{\corpos \symbol{80}}}        
\newcommand{\R}{\mbox{\corpos \symbol{82}}}          
\newcommand{\T}{\mbox{\corpos \symbol{84}}}          

\newcommand{\g}{\mbox{\lie g}}       %

\newcommand{\rr}{\rightarrow}
\newcommand{\lrr}{\longrightarrow}

\newcommand{\hnab}{{\cal H}^{\nabla}}             %

\newcommand{\Tr}[1]{{\mathrm{Tr}}\,{#1}}

\newcommand{\End}[1]{{\mathrm{End}}\,{#1}}
\newcommand{\Ad}[1]{{\mathrm{Ad}}\,{#1}}

\newcommand{\Id}{{\mathrm{Id}}}
\newcommand{\dx}{{\mathrm{d}}}
\newcommand{\inv}[1]{{#1}^{-1}}

\newcommand{\cinf}[1]{{\mathrm{C}}^\infty_{#1}}

\newcommand{\vol}{{\mathrm{Vol}}\,}
\newcommand{\ric}{{\mathit{Ric}}\,}

\newcommand{\nastar}{\nabla^\star}

\newtheorem{teo}{Theorem}[section]

\newtheorem{coro}{Corollary}[section]
\newtheorem{prop}{Proposition}[section]

\numberwithin{equation}{section}

\def\cyclic{\mathop{\kern0.9ex{{+}
\kern-2.2ex\raise-.28ex\hbox{\Large\hbox{$\circlearrowright$}}}}\limits}

\begin{document}
\thispagestyle{empty}

\begin{center}

{\huge{\bf{The $G_2$ sphere over a 4-manifold}}}

\vspace{1cm}

R. Albuquerque\footnote{Departamento de Matem\'atica da Universidade de \'Evora and Centro de Investiga\c c\~ao em Matem\'atica e Aplica\c c\~oes (CIMA), Rua Rom\~ao Ramalho, 59, 7000 \'Evora, Portugal.}
\hspace{4.85cm}
I. M. C. Salavessa\footnote{Centro de F\'\i sica das Interac\c c\~oes Fundamentais (CFIF), Instituto Superior T\'ecnico, Av. Rovisco Pais, Edif\'\i cio Ci\^encia, Piso 3, 1049-001 Lisboa, Portugal.}

rpa@uevora.pt\hspace{4cm}  isabel.salavessa@ist.utl.pt

\vspace{1.3cm}

\date{\today} 

\begin{small}{\bf Abstract}
\end{small}
\end{center}

\begin{small}

We present a construction of a canonical $G_2$ structure on the unit sphere tangent bundle $S_M$ of any given orientable Riemannian 4-manifold $M$. Such structure is never geometric or 1-flat, but seems full of other possibilities. We start by the study of the most basic properties of our construction. The structure is co-calibrated if, and only if, $M$ is an Einstein manifold. The fibres are always associative. In fact, the associated 3-form $\phi$ results from a linear combination of three other {\it volume} 3-forms, one of which is the volume of the fibres. We also give new examples of co-calibrated structures on well known spaces. We hope this contributes both to the knowledge of special geometries and to the study of 4-manifolds.
\end{small}

\vspace*{4mm}


\vspace*{2mm}


\vspace*{10mm}

The authors acknowledge the support of Funda\c{c}\~{a}o Ci\^{e}ncia e Tecnologia, either through the project POCI/MAT/60671/2004 and through their research centers, respectively, CIMA and CFIF.

\markright{\sl\hfill  Albuquerque --- Salavessa \hfill}

\vspace*{10mm}

\section{Introduction}
\label{Introduction}

Since M. Berger found the famous list of possible holonomy groups for locally irreducible Riemannian metrics, the exceptional Lie group $G_2$ has won a new importance to geometry and physics and the quest for manifolds with such holonomy has become a very defying challenge. Nevertheless, it took some time until the first dedicated articles appeared (\cite{Bonan}) on $G_2$ manifolds, ie. 7 dimensional manifolds having a Riemannian metric with holonomy group inside $\mathrm{Aut}(\Octoni)$.

Since $G_2=\{g\in SO(7)|\ g^*\phi=\phi\}$, where $\phi$ is a non-degenerate 3-form on $\R^7$, the structure is defined by a non-degenerate smooth differential 3-form on the Riemannian manifold. Reciprocally, $\phi$ and the orientation determine the metric. When this structure is `geometric' or 1-flat, ie. when $\phi$ is parallel for the Levi-Civita connection, the manifold is called a $G_2$ manifold. The first examples with complete metrics of exactly $G_2$ holonomy were found by R. Bryant on open subsets of $\R^7$ (\cite{Bryant1}). Soon afterwords, with S. Salamon, in \cite{BrySal}, they applied a new set of ideas to other spaces with non-degenerate 3-forms. For instance, on the product of a 3-sphere with $\R^4$ and on the vector bundle of self dual 2-forms over an Einstein $4$-manifold. The first compact examples with holonomy the whole $G_2$ were constructed by D. Joyce, \cite{Joy}, using mainstream techniques from calculus on manifolds.

A link to modern unifying theories of physics is refered in works of M. F. Atiyah and E. Witten, \cite{AtiWit}, and Th. Friedrich and S. Ivanov, \cite{FriIva1,FriIva2}. The former article contains various examples of fibre bundles carrying a metric with holonomy inside $G_2$. Namely rank $n$ vector bundles over $(7-n)$-manifolds, but also those other examples of developing cone singularities from $SU(3)$ six-manifolds, in particular, the twistor spaces of $S^4$ and $\C\Pro^2$, which have already produced true $G_2$ manifolds. Important to recall too are the $S^1$ products with Calabi-Yau spaces. Or the solvable $\R$ extensions built from 6-nilmanifolds, by I. Agricola, S. Chiossi, A. Fino, S. Salamon and others, cf. \cite{AgriChiFino,ChiFino}. Specially interesting to refer here is the article \cite{BobNur} with respect to constructions of $G_2$ structures on sphere bundles. They present their ``twistor'' non-degenerate 3-forms defined on $S^2$-bundles over Riemannian 5-manifolds which admit a certain special geometry, a reduction to $SO(3)$.

Most $G_2$ structures known today are not `flat to first order'. (As explained in \cite{Bryant2}, recall that a local correspondence by diffeomorphism with the \textit{infinitly flat} model $\R^7$, leading to a complete integrability, is a much more difficult question). The interest of co-calibrated structures, defined by the weaker condition $\delta\phi=0$, is proved for instance in recent work of Th. Friedrich and S. Ivanov, \cite{FriIva2}. The stricter case of nearly parallel structures, $\dx\phi=c*\phi$, with $c$ constant, has also been thoroughly studied, as we may see in \cite{FriKaMoSe}.

A manifold admits a $G_2$ structure if, and only if, its two first Stiefel-Whitney classes $w_1,w_2$ vanish. The present article deals with a new and particular $G_2$ structure corresponding with any given orientable Riemannian 4-manifold $M$; no further assumptions being made. The space is the $S^3$-bundle of unit tangent vectors $S_M\subset TM$, which inherits the usual orientation and metric from an orthogonal decomposition of the tangent space of $M$ into horizontal and vertical vector bundles. The question of reduction to $G_2$ is then solved. We prove we cannot have such strict condition and analyse the associated co-calibration, which happens if, and only if, $M$ is an Einstein manifold.

As the reader may care to notice, our construction is reminiscent of the twistor definition of an almost complex structure on the bundle of linear complex structures of even dimensional manifolds. (Hence the purpose of the introduction above, with such punctual references to the history of $G_2$'s.) Not only we apply here some of the techniques of the well known theory of twistors, but also the structure can be seen as one of those following from a tautological definition. There are even good expectations that the ``$G_2$ sphere" introduced in this article beholds other parallels with the celebrated Penrose-Ward correspondence for the study of 4 dimensional manifolds.

We acknowledge Professor D. V. Alekseevsky for fruitful conversations and for calling our attention to the importance of $G_2$ geometry.

\section{Riemannian geometry of the sphere bundle}
\label{rrg}

\subsection{Recalling the pull-back connection}
\label{rtpbc}

We start by recalling some material from the theory of connections. Suppose we have a smooth vector bundle $\pi:E\rr M$ over a smooth manifold $M$ together with a connection $\nabla$ on $E$. Suppose also that it is associated to a principal $G$-bundle $F$, where $G$ is a Lie group, so that we may write $E=F\times_GV$, with $V$ a space endowed with a left $G$-module structure. The connection is given by a 1-form $\alpha$ on $F$ taking values on the Lie algebra of $G$ and such that
\begin{equation}\label{eqdeconexao1}
  R_g^*\alpha=\Ad{(\inv{g})}\alpha,\hspace{2cm}\alpha_p\Bigl(\frac{\dx}{\dx t}_{|t=0}\,p\exp tA\Bigr)=A ,
\end{equation}
$\forall g\in G,\ A\in\mbox{Lie}(G)$. Any section of $E$ may be written locally, say on some open subset $U$ of $M$, as a product $sv$ with $s\in\Gamma(U;F):=\Omega_U(F),\,\ v\in \cinf{U}(V)$. Then we may state the formula relating the above:
\begin{equation}\label{eqdeconexao2}
 \nabla_X \,sv=s(s^*\alpha)(X)v+s\dx v(X),
\end{equation}
$\forall X\in TM$. Recall that the kernel of $\alpha$ is a
tangent distribuition on $F$, usually called the space of horizontal
directions, isomorphic to $\inv{\pi}TM$.

Now we assume further that we are given a map $f:Z\rr M$ through which the submersion $\pi:F\rr M$ factors, ie. there exists another submersion $\pi_1$ such that $f\circ\pi_1=\pi$. Then we have the following result, on whose proof we cast some light.
\begin{prop}\label{propo1}
The vector bundle $f^*E\rr Z$ is associated to $\pi_1:F\rr Z$ as a principal $H$-bundle, where $H=\{g\in G:\ \pi_1(pg)=\pi_1(p),\ \forall p\in F\}$. Moreover,
\begin{equation}\label{eq3}
(f^*\nabla)\,\xi v=\xi(\xi^*\alpha)\:v+\xi\dx v
\end{equation}
for any section $\xi v$, where $\xi\in\Gamma({\cal U};F),\ v\in\cinf{{\cal U}}(V),\ {\cal U}$ open in $Z$.
\end{prop}
The result shows how a reduction to the structure group $H$ spontaneously occurs. It follows from the commutative diagram
\begin{equation*}\label{eq4}
\begin{array}{ccccc}
F & \stackrel{\epsilon}\lrr & f^*F &\stackrel{\mbox{\tiny{$\mathrm{pr}_2$}}}\lrr &F \\
\pi_1\downarrow &        &\downarrow\mbox{\tiny{$\mathrm{pr}_1$}} & & \downarrow \\
Z              &  \lrr   & Z& \lrr& M
\end{array}
\end{equation*}
where $\epsilon(p)=(\pi_1(p),p)$, clearly an $H$-equivariant map. Formula (\ref{eq3}) follows from standard computations evaluating the new connection form, supported on relations (\ref{eqdeconexao1},\ref{eqdeconexao2}); more easily understood if we see $\xi=f^*s\,g=(s\circ f)g$ for some $s\in\Gamma(\pi_1({\cal U});F),\ g\in\cinf{{\cal U}}(G)$. It is interesting to notice that if $Y\in TZ$ is vertical, ie. $\dx f(Y)=0$, then \begin{equation}\label{eq5}
(f^*\nabla)_Yf^*s\,g=f^*s\,\dx g(Y).
\end{equation}
Finally, we recall how the curvature tensors relate: $R^{f^*\nabla}=f^*R^{\nabla}$.

\subsection{The sphere bundle and its Levi-Civita connection}
\label{tsbailcc}

Let $(M,\langle\:,\:\rangle)$ be a smooth orientable Riemannian manifold of dimension $m$. Let
\begin{equation}\label{spherebundle} 
S_M=\bigl\{u\in TM|\ \|u\|=1\bigr\},
\end{equation}
which is the total space of a sphere bundle over $M$, and let $f$ denote the projection to the base. Notice $S_M$ is always orientable. There is a well known exact sequence $0\rr {\cal V}\rr T\,S_M\rr f^*TM\rr 0$ lying above $S_M$ and the Levi-Civita connection of $M$ induces a direct sum decomposition $T\,S_M={\cal V}\oplus\hnab$ into a rank $m$ horizontal distribuition and a vertical ${\cal V}=\ker\dx
f$ distribuition. While the latter is integrable because it agrees with the tangent to the fibres, hence is closed under Lie bracket of vector fields, the former will be integrable depending on the vanishing of the Riemann curvature tensor of $M$.

Notice the bundle $\cal V$ may be identified with the subvector bundle of $\inv{f}TM$ (we use the notation $\inv{f}$ to refer to the vertical side) such that ${\cal V}_u=u^\perp\subset T_{f(u)}M$. We denote by $U$ the ubiquous section of the pullback bundle such that $U_u=u$. Then we may also write $U^\perp$ for the vertical bundle. Virtually, $U$ appears both on the vertical and the horizontal sides of $TTM$ but we shall see it as an object in the vertical side.
\begin{prop}\label{propo2}
$\hnab=\{X\in TS_M|\ (f^*\nabla)_XU=0\}$ and
$f^*\nabla_{\,Y}U=Y$ for any $Y\in{\cal V}$.
\end{prop}
\begin{proof}
Consider the orthonormal frame bundle
$SO_M=\{p:\R^m\rr T_xM\,|\ p\ \mbox{isometry}, \ x\in M\}$ of $M$, consider a vector $v_0\in S^{m-1}$ and a fibre bundle $\pi_1:SO_M\rr S_M$ defined by $\pi_1(p)=pv_0$. Let $u_0=p_0v_0\in S_M$. Then $\inv{f}TM$ is associated to $SO_M$ as a principal bundle and thus $U=\xi v_0$ for some local section $\xi$. By proposition \ref{propo1} we know there is a horizontal $\hnab\subset T_{u_0}S_M$ coming from the $SO(m-1)$-invariant distribuition $\ker\alpha$.

Now suppose $X_{u_0}={\pi_1}_*\tilde{X}$ is horizontal, ie. with $\tilde{X}\in\ker\alpha_{p_0}$. It is clear that the horizontal part of $\xi_*(X_{u_0})$ is $\tilde{X}$ because ${\pi_1}_*$ is a bijection on the horizontals. But by deeper results in connection theory (though one may argue just with integral curves), we may further assume that we have a section $\xi$ on a neighborhood of $u_0$ such that $\xi(u_0)=p_0$ and $\xi_*(X_{u_0})=\tilde{X}_{p_0}$. Therefore $f^*\nabla_{X}\,U=\xi(\xi^*\alpha)(X)v_0+\xi\dx v_0(X)=0$ by equation (\ref{eq3}). Suppose now $Y\in\ker\dx f_{u_0}$ is vertical; then we find by equation (\ref{eq5})
\[ f^*\nabla_Y\:U=\xi\,\dx g(Y)v_0=\dx(s_{f(u_0)}gv_0)(Y)=\dx(u\mapsto u)_{u_0}(Y)=Y \]
since we may restrict to the fibre and write $\xi=f^*s\,g$ with $s$ a local section of
$SO_M\rr M$ and $g\in C^\infty(SO(m))$.
\end{proof}
The previous result is quite used in the literature, though oftenly not proved: it is well established that the vertical part of $\dx X(Z)$ is precisely $\nabla_ZX$, for a given section $X\in\Omega(TM)$. Indeed, in this case $(\dx X(Z))^v=f^*\nabla_{\dx X(Z)}U=(X^*f^*\nabla)_ZX^*U=\nabla_ZX$.

Now, we may endow $S_M$ with a Riemaniann structure. Naturally, it is given by $f^*\langle\:,\:\rangle$ on $\hnab$ identified isometrically with the pull-back bundle $f^*TM$  and it is $f^*\langle\:,\:\rangle$, again, on the restriction to $\cal V$. Surely the two distribuitions sit orthogonally inside $TS_M$.

The tangent bundle clearly inherits a metric connection preserving the decomposition ${\cal V}\oplus\hnab$, which we shall denote by $\nastar$. Notice on the vertical tangent directions we must add a correction term
\begin{equation}\label{correcterm}
\nastar X^v=f^*\nabla X^v-\langle f^*\nabla X^v,U\rangle U
\end{equation}
as it is well known. Let $R^*=f^*R^{\nabla}$ denote the curvature of $f^*\nabla$. We see $R^*U\in\Omega^2({\cal V})$.

We have used $\,\cdot\,^v:TS_M\rr{\cal V}$ to denote the projection,
$X^v=(f^*\nabla)_XU=\nastar_XU$. The notation $X^h$ corresponds to the horizontal
part of the vector, which ocasionally is identified with $\dx f(X)$ too.
\begin{teo}
The Levi-Civita connection of $S_M$ is given by
\begin{equation}
D_XY=\nastar_XY-\dfrac{1}{2}R^*_{X,Y}U+{\cal A}_XY
\end{equation}
where ${\cal A}$ takes values in $\hnab$ and is defined by
\begin{equation}
\langle {\cal A}_XY,Z^h\rangle=\dfrac{1}{2}\langle R^*_{X^h,Z^h}U, Y^v\rangle + \dfrac{1}{2}\langle R^*_{Y^h,Z^h}U, X^v\rangle.
\end{equation}
\end{teo}
\begin{proof}
Let us first see the horizontal part of the torsion:
\[ \dx f(T^D(X,Y)) = \nastar_XY^h+{\cal A}_XY-\nastar_YX^h-{\cal A}_YX-\dx f[X,Y] = f^*T^\nabla(X,Y)\ =\ 0   \]
because this is how the torsion tensor of $M$ lifts to $f^*TM$ and since ${\cal A}$ is symmetric. Now we check the vertical part.
\begin{eqnarray*}
(T^D(X,Y))^v& = &\nastar_XY^v-\dfrac{1}{2}R^*_{X,Y}U-\nastar_YX^v+ \dfrac{1}{2}R^*_{Y,X}U-[X,Y]^v\\
& = &\nastar_X\nastar_YU-R^*_{X,Y}U-\nastar_Y\nastar_XU-\nastar_{[X,Y]}U\ =\ 0.
\end{eqnarray*}
It remains to check $D$ is a metric connection, which is equivalent
to the difference with $\nastar$ being skew-adjoint. This is an easy straightforward computation, so we are finished.
\end{proof}
We remark that $D\hnab\subset\hnab$ if, and only if, the Riemmannian manifold $M$ is flat. We also note that $D$ is the Levi-Civita connection of $TM$ up to the correction term referred in (\ref{correcterm}).

\section{The canonical $G_2$ structure}
\label{ag2s}

\subsection{The octonionic line}
\label{alos}

Let $Q$ be a 4-dimensional oriented Euclidian vector space and let $u$ denote a fixed vector in $Q$ with norm 1. This is sufficient to define a unique quaternionic structure on $Q$. Indeed, we have a canonical vector cross product on $u^\bot$ given by $\langle X\times Y,Z\rangle=\vol(u,X,Y,Z)$ and hence the quaternionic product:
\begin{equation}\label{quatpro}
(\lambda u+X)(\mu u+Y)=(\lambda\mu-\langle X,Y\rangle)u+\lambda Y+\mu X+X\times Y
\end{equation}
for any $\lambda,\mu\in\R,\ X,Y\in u^\perp$. Having a conjugation map $\overline{\lambda u+X}=\lambda u-X$, we may proceed to establish an octonionic structure in $T=Q\oplus Q$:
\begin{equation}\label{octoprod}
(a_1,a_2)(a_3,a_4)=(a_1a_3-\overline{a_4}a_2,a_4a_1+a_2\overline{a_3})
\end{equation}
for all $a_1,\ldots,a_4\in Q$. This is the well known Cayley-Dickson process (cf. \cite{HarLaw}).
Recall the cross product of two imaginary quaternions $X,Y\in u^\bot$ is equal to the imaginary part of $\overline{Y}X$ and notice the following formula for the product of two imaginary octonions. If $a_i=\lambda_iu +X_i,\ i=2,4$, and $X_1,\ldots,X_4\in u^\bot$, then
\begin{equation}\label{imagoctoprod}
\begin{split}
(X_1,a_2)(X_3,a_4)\  =\ (X_1X_3-\overline{a_4}a_2,a_4X_1-a_2X_3) \ = \hspace{3cm}\\
 \  = \ \bigl( (-\lambda_4\lambda_2- \langle X_4,X_2\rangle-\langle X_1,X_3\rangle)u+X_1\times X_3-\lambda_4 X_2+\lambda_2 X_4+X_4\times X_2, \bigl. \\
 \ \ \ \ \ \bigr. (\langle X_2,X_3\rangle-\langle X_4,X_1\rangle)u +\lambda_4X_1+X_4\times X_1-\lambda_2 X_3-X_2\times X_3\bigr).\ \ \ 
\end{split}
\end{equation}
Finally we get a non-degenerate 3-form in $u^\bot\oplus Q\subset T$:
\begin{equation}\label{formalinear}
\begin{split}
\lefteqn{\phi((X_1,a_2),(X_3,a_4),(X_5,a_6))\  =\ \langle(X_1,a_2)(X_3,a_4),(X_5,a_6)\rangle \ = \ \ \ }\\
& \ \  = \ \langle X_1\times X_3,X_5\rangle-\lambda_4 \langle X_2,X_5\rangle  +\lambda_2\langle X_4,X_5\rangle+\langle X_4\times X_2,X_5 \rangle
+\lambda_6\langle X_2,X_3\rangle \\
& \ \ \ -\lambda_6\langle X_4,X_1\rangle+\lambda_4\langle X_1,X_6 \rangle
+\langle X_4\times X_1,X_6\rangle -\lambda_2\langle X_3,X_6\rangle
-\langle X_2\times X_3,X_6\rangle.
\end{split}
\end{equation}
It is also important to recall that $X\times$ is a skew-symmetric operator in $u^\bot$.

\subsection{The $G_2$ structure on the sphere bundle}
\label{tg2sotsb}

Now let us go back to the setting where we have a Riemannian manifold $M$, assume it is 4-dimensional and orientable, and consider its sphere bundle. Using the Levi-Civita connection to produce a splitting of the tangent bundle of $S_M$ and identifications with $\inv{f}TM$, the previous construction of a linear octonionic structure yields a $G_2=\mathrm{Aut}(\Octoni)$ structure on $TS_M$. The canonical section $U$ of $\inv{f}TM$ plays the role of the real part in the normed ring.

In order to study the structure thus presented, we shall need a copy of $U$ in $\hnab$ as well as the metric. Therefore it is wise to introduce the isomorphism
\begin{equation}\label{teta}
\theta:\hnab\lrr\inv{f}TM\supset{\cal V}
\end{equation}
defined by $f^*\nabla$-parallel and isometric identifications with the pull-back bundle of $TM$, explained in section \ref{tsbailcc}. Since we want to differentiate under the Levi-Civita connection $D$, we extend $\theta$ by 0 to the vertical tangent bundle, taking values in the octonionic vector bundle on $S_M$. So in fact we have
\begin{equation}\label{teta1}
\theta:TS_M\rr TTM
\end{equation}
as a kind of a ``soldering form".

A $G_2$ structure is entirely determined by a non-degenerate 3-form. We shall be using a few smooth scalar tensors over $S_M$ in order to determine and study the 3-form of the present example, say $\phi$. In place of the cross product, we define $\alpha=U\lrcorner\,\inv{f}\vol_M\in\Omega^0(\Lambda^3{\cal V}^*)\subset\Omega^3$. We also need a 1-form $\mu$ defined by $\mu(X)=\langle U,\theta(X)\rangle$ and a 2-form $\beta$ such that $\beta(X,Y)=\langle\theta X,Y\rangle-\langle\theta Y,X\rangle=\langle\theta X^h,Y^v\rangle-\langle\theta Y^h,X^v\rangle$. 

Let us also establish some notation for not so well established computations. Given any $p$-tensor $\eta\in\otimes^pT^*M$ and any endomorphisms $B_i$ of the tangent bundle we let $\eta\circ (B_1\wedge\ldots\wedge B_p)$ denote the new $p$-tensor defined by
\begin{equation}\label{esttensorcont}
\eta\circ (B_1\wedge\ldots\wedge B_p)(Y_1,\ldots,Y_p)=\sum_{\sigma\in S_p}\mathrm{sg}(\sigma)\eta(B_1Y_{\sigma_1},\ldots,B_pY_{\sigma_p}).
\end{equation}
It is easy to check such contraction is parallel and thus that it obeys a simple Leibniz rule under covariant differentiation, with no $-1$ signs attached. For instance, if $\eta$ is a $p$-form, then $\eta\circ\wedge^p\Id=p!\,\eta$. Furthermore, one verifies that for a wedge of 1-forms, $\eta_1\wedge\ldots\wedge\eta_p\circ(B_1\wedge\ldots\wedge B_p)= \sum_{\sigma\in S_p} \eta_1\circ B_{\sigma_1}\wedge\ldots\wedge\eta_p\circ B_{\sigma_p}$.

Using the above, we define $\alpha_i\in\Omega^3$, for $i=1,2$, by
\begin{equation}\label{alfas} 
\alpha_1=\frac{1}{2}\alpha\circ(\theta\wedge\Id\wedge\Id),\hspace{1.7cm}
\alpha_2=\frac{1}{2}\alpha\circ(\theta\wedge\theta\wedge\Id).
\end{equation} 
Finally the associated 3-form $\phi$ of the $G_2$ structure on the sphere bundle of $M$ induced from (\ref{formalinear}) satisfies
\begin{equation}\label{threeform}
\begin{split}
\phi(X,Y,Z)\  = \ \langle X^v\times Y^v,Z^v\rangle- \langle Y^h,U\rangle\langle X^h,Z^v\rangle  +\langle X^h,U\rangle\langle Y^h,Z^v\rangle+  \\
 \langle Y^h\times X^h,Z^v \rangle +\langle Z^h,U\rangle\langle X^h,Y^v\rangle- \langle Z^h,U\rangle\langle Y^h,X^v\rangle +  \\
\langle Y^h,U\rangle\langle Z^h,X^v\rangle -\langle Y^h\times Z^h,X^v\rangle -\langle X^h,U\rangle\langle Z^h,Y^v\rangle +\langle X^h\times Z^h,Y^v\rangle,
\end{split}
\end{equation}
written in a heuristic if not confusing way. Since in fact $\langle X^h,Z^v\rangle$ is given by $\langle\theta X,Z\rangle$ in the new framework, we may claim to have found the simple expression
\begin{equation}\label{threeformagain}
\phi=\alpha+\mu\wedge\beta-\alpha_2.
\end{equation}
\begin{prop}[basic structure equations] \label{bse}
We have the following basic relations:
\begin{equation}\label{bse1} 
\begin{split}
*\alpha=f^*\vol_M,\ \ \ \ \ *\alpha_1=-\mu\wedge\alpha_2,\ \ \ \ \ *\alpha_2=\mu\wedge\alpha_1,\hspace{2.4cm} \\
*\beta=-\frac{1}{2}\mu\wedge\beta^2,\ \ \ \ \ *\beta^2=-2\mu\wedge\beta,\ \ \ \ \ \beta^3\wedge\mu=-6\vol_{S_M},\hspace{2.4cm} \\
\alpha_1\wedge\alpha_2=3*\mu=-\frac{1}{2}\beta^3,\ \ \ \ \ \beta\wedge\alpha_i=\beta\wedge*\alpha_i=\alpha_0\wedge\alpha_i=0,\ \ \forall i=0,1,2,
\end{split}
\end{equation}
where we denote $\alpha=\alpha_0$. Henceforth $\alpha\wedge\phi=\alpha_2\wedge\phi =*\alpha_1\wedge\phi=0$. Moreover, $*\phi=f^*\vol_M-\frac{1}{2}\beta^2-\mu\wedge\alpha_1$ and thus $*\alpha\wedge\phi=\alpha\wedge*\phi=\vol_{S_M}$.
\end{prop}
\begin{proof}
An easy way to study all these forms is by reference to a frame.
Let $e_0=\theta^tU,e_1,\ldots,e_6$ denote a direct orthonormal basis of $TS_M$. So we assume $e_0,\ldots,e_3$ is a direct orthonormal basis of $\hnab$ and $e_{i+3}=\theta e_i, \ \forall i=1,2,3$. By definition, $\alpha=e^{456}$ and $\mu=e^0$. Thus $*\alpha=e^{0123}=f^*\vol_M$. It is also trivial to see $\beta=e^{14}+e^{25}+e^{36}$. From direct inspection on $\alpha$ composed with $\theta$s  we easily find
\begin{equation}\label{alphascoef}
  \alpha_1=e^{156}+e^{264}+e^{345}\ \ \ \ \mbox{and}\ \ \ \ \alpha_2=e^{126}+e^{234}+e^{315} .
\end{equation}
Hence $*\alpha_1=-e^{0234}+e^{0135}-e^{0126}=-\mu\wedge\alpha_2$ and $*\alpha_2=e^{0345}+e^{0156}+e^{0264}=\mu\wedge\alpha_1$. Now
\[ \beta^3=(2e^{1425}+2e^{1436}+2e^{2536})\wedge\beta=6e^{142536}=-6*\mu . \]
Finally, $*\beta=e^{02356}+e^{01346}+e^{01245}=-\frac{1}{2}\mu\wedge\beta^2$ and  $*\beta^2=-2*(e^{1245}+e^{1346}+e^{2356})=-2(e^{036}+e^{025}+e^{014})$ and so the result follows. To understand the last relations in (\ref{bse1}) one just has to look to the formulae written along this proof.
\end{proof}
The reader may check directly $|\phi|^2=7$. Now we need the computation of some derivatives.
\begin{prop}\label{derivacoes}
For any vector field $X$ over $S_M$:\\
1. $D_X\alpha=\frac{1}{4}\alpha\circ(R^*_{X,\cdot}U\wedge\Id\wedge\Id)={\cal A}_X\alpha$.\\
2. $2D_X\alpha_2=D_X(\alpha\circ\theta\wedge\theta\wedge\Id)=({\cal A}_X\alpha)\circ\theta\wedge\theta\wedge\Id-2\alpha\circ\,\theta {\cal A}_X\wedge\theta\wedge\Id$.\\
3. $D_X\mu=X^\flat\circ\theta -\mu\circ {\cal A}_X$.\\
4. $\dx\mu=-\beta$ and $\delta\mu=0$.
\end{prop}
\begin{proof}
1. We have $D_XY_i=\nastar_XY_i-\frac{1}{2}R^*_{X,Y_i}U+{\cal A}_XY_i$ for
any three vector fields $Y_1,Y_2,Y_3$ on $S_M$ and thus
\begin{eqnarray*}
D_X\alpha(Y_1,Y_2,Y_3)\ =\ (\nastar_X\inv{f}\vol_M)(U,Y_1,Y_2,Y_3)
+\inv{f}\vol_M(\nastar_XU,Y_1,Y_2,Y_3)  \\
 +\frac{1}{2}\bigl(\alpha(R^*_{X,Y_1}U,Y_2,Y_3) +\alpha(Y_1,R^*_{X,Y_2}U,Y_3) +\alpha(Y_1,Y_2,R^*_{X,Y_3}U)\bigr).
\end{eqnarray*}
The first two terms on the sum vanish because $\nabla\vol_M=0$ and because $\nastar_XU$ and the $Y_i^v\in\cal V$ are linearly dependent. Hence the result. If we see $\alpha=Y^{123}=Y^1\wedge Y^2\wedge Y^3$, where the $Y_i$ form an orthonormal basis of $\cal V$, then 
\[ \frac{1}{4}\alpha\circ(R^*_{X,\cdot}U\wedge\Id\wedge\Id)=
({\cal A}_XY_1)^\flat\wedge Y^{23}-({\cal A}_XY_2)^\flat\wedge Y^{13}+({\cal A}_XY_3)^\flat\wedge Y^{12}, \]
ie. ${\cal A}_X$ acts as a derivation of $\alpha$. \\
2. Let $\tilde{D}=D+\langle f^*\nabla\ ,U\rangle U$, the
Levi-Civita connection of $TM$. Since $\alpha$ is 0 when we take one
direction proportional to $U$, the $D$-derivative we have to compute
can be done with $\tilde{D}$. This gives us the possibility of
computing $\tilde{D}_X\theta$, which results in
$[{\cal A}_X,\theta]$. Hence
\begin{eqnarray*}
D_X(\alpha\circ\theta\wedge\theta\wedge\Id) & =& D_X\alpha\circ\, \theta\wedge\theta\wedge\Id +
  2\alpha\circ\,\tilde{D}_X\theta\wedge\theta\wedge\Id \\
&=& ({\cal A}_X\alpha)\circ\theta\wedge\theta\wedge\Id-2\alpha\circ\,
\theta {\cal A}_X\wedge\theta\wedge\Id.
\end{eqnarray*}
3. Here we just have to compute: $(D_X\mu)Y=X(\mu Y)-\mu(D_XY)=\langle f^*\nabla_XU,\theta Y\rangle +\langle U,f^*\nabla_X(\theta Y)\rangle-\langle U,\theta(f^*\nabla_X Y)\rangle-\langle U,\theta({\cal A}_XY)\rangle=\langle X,\theta Y\rangle-\mu({\cal A}_XY)$. \\
4. $\dx\mu=-\beta$ is a simple computation arising from $\dx\mu(X,Y)=(D_X\mu)Y-(D_Y\mu)X$ and from the symmetry of $\cal A$. Furthermore, $\delta\mu=-*\dx *\mu=\frac{1}{6}*\dx\beta^3=0$.
\end{proof}
Of course, a formula for $D\alpha_1$ follows as in 2 above.

We recall that an orientable Riemannian 7-manifold with a $G_2$ structure $\phi$ admits a (holonomy) reduction to a subgroup lying in $G_2$ if, and only if, $\phi$ is parallel. Such condition being fullfield gives place to the concept of a $G_2$-manifold. Furthermore, by a result of \cite{Fer}, this is equivalent to having $\phi$ harmonic. If $\dx\phi=0$, then the structure is called calibrated and, if $\delta\phi=0$, the structure is known as co-calibrated.

Recall the Ricci tensor of $M$ is defined by $r(X,Y)=\Tr{R^\nabla_{\cdot,X}Y}$. It is also given by a symmetric endomorphism $\ric\in\Omega(\End{TM})$ satisfying $r(X,Y)=\langle \ric X,Y\rangle$. These tensors lift to $S_M$ in the usual way, through $f^*$ or $\inv{f}$. We shall see $(\ric U)^\flat\in\Omega({\cal V}^*)$ as a 1-form vanishing on $\hnab$ and restricted to vertical tangent directions.
\begin{teo}\label{dphiemframe}
With a frame $e_0,\ldots,e_6$ such that $e_0=\theta^tU,\ e_{i+3}=\theta e_i$, $i=1,2,3$, induced from an oriented orthonormal frame of $M$, and setting ${\cal R}^{ij}=\langle R^\nabla_{\cdot,\cdot}e_i,e_j\rangle=\sum_{0\leq k<l\leq3}{\cal R}^{ij}_{kl}e^{kl}$, we have
\begin{equation}
\dx\phi\ =\ {\cal R}^{01}\wedge e^{56}+{\cal R}^{02}\wedge e^{64}+{\cal R}^{03}\wedge e^{45}-\beta^2+r(U,U)f^*\vol_M.
\end{equation}
\end{teo}
\begin{proof}
Since $\dx\phi=\dx(\alpha-\alpha_2)-\beta^2$, we start by looking at $\dx\alpha_2$. By proposition \ref{derivacoes} we find
\[ \dx\alpha_2=\sum_{i=0}^6e^i\wedge D_{e_i}\alpha_2=\frac{1}{2}\sum_{i=0}^6e^i\wedge \bigl( ({\cal A}_{e_i}\alpha)\circ\theta\wedge\theta\wedge\Id -2\alpha\circ\,\theta {\cal A}_{e_i}\wedge\theta\wedge\Id\bigr) . \]
In the second term on the right we have a contraction of a symmetric derivation $\cal A$ within a skew tensor, so it is easy to see that it vanishes. Since $\alpha=e^{456}$ and since $\cal A$ takes only horizontal values, we have
\begin{eqnarray*}
\dx\phi &=&\sum_{i=0}^3\bigl(e^i\wedge({\cal A}_{e_i}e_4)^\flat\wedge e^{56}+ e^i\wedge({\cal A}_{e_i}e_5)^\flat\wedge e^{64}+ e^i\wedge({\cal A}_{e_i}e_6)^\flat\wedge e^{45}  \\
& & \bigl. -e^i\wedge({\cal A}_{e_i}e_4)^\flat\wedge e^{23} 
-e^i\wedge({\cal A}_{e_i}e_5)^\flat\wedge e^{31}-e^i\wedge({\cal A}_{e_i}e_6)^\flat\wedge e^{12}\bigr)- \beta^2.
\end{eqnarray*}
We know that for $i,k\leq3<j$ we have by definition $a_{ijk}=\langle{\cal A}_{e_i}e_j,e_k\rangle=\tfrac{1}{2}\langle R^\nabla_{e_i,e_k}e_0,e_{j-3}\rangle$. Thus $a_{ijk}=-a_{kji}$ and we deduce $\dx\phi$ also takes the shape
\begin{eqnarray*}
 2\sum_{0\leq i<j\leq 3}\bigl((a_{i4j}e^{ij}+\tfrac{1}{6}e^{23})\wedge e^{56}+(a_{i5j}e^{ij}+\tfrac{1}{6}e^{31})\wedge e^{64}+(a_{i6j}e^{ij}+\tfrac{1}{6}e^{12})\wedge e^{45}\bigr) \\
 -2(a_{041}+a_{052}+a_{063})e^{0123}.  \hspace{3.5cm}
\end{eqnarray*}
That the last term is essentially the Ricci curvature of $M$ is easily checked. 
\end{proof}
\begin{coro}\label{ng2man}
$(S_M,\phi)$ is never a $G_2$-manifold.
\end{coro}
\begin{proof}
If we compute $\dx\phi_{0156}$ in $u\in S_M$, we find $2a_{041}=\langle R^\nabla_{u,e_1}u,e_1\rangle$. Since $u,e_1$ may be chosen such to span any plane in $TM$ we conclude that the assumed calibration of $S_M$ implies $M$ flat. However, this yields $\dx\phi=-\beta^2$ in contradiction with the hypothesis.
\end{proof}
\begin{prop}\label{cocali}
We have $\dx*\phi=-f^*\vol_M\wedge(\ric U)^\flat$. In particular,
$S_M$ is co-calibrated if, and only if, $M$ is an Einstein manifold.
\end{prop}
\begin{proof}
From propositions \ref{bse},\ref{derivacoes} we immediately find $\dx *\phi=-\dx\mu\wedge\alpha_1+\mu\wedge\dx\alpha_1=\mu\wedge\dx\alpha_1$.
As in the computation of $\dx\alpha_2$, we now get $\dx*\phi=$
\begin{eqnarray*}
&=& \sum_{i=0}^6\mu\wedge e^i\wedge\, \tfrac{1}{2}
({\cal A}_{e_i}\alpha)\circ\theta\wedge\Id\wedge\Id \ =\ \sum_{i=0}^3\mu\wedge 
 e^{i}\wedge\bigl(({\cal A}_{e_i}e_4)^\flat\wedge(e^{26}+e^{53}) \\
 & &\ \ \ -({\cal A}_{e_i}e_5)^\flat\wedge(e^{16}+e^{43}) +({\cal A}_{e_i}e_6)^\flat\wedge(e^{15}+e^{42}) \bigr)\\
&=& 2\mu\wedge\sum_{1\leq i<j\leq 3}\bigl((a_{i4j}(e^{ij26}+e^{ij53}) +a_{i5j}(e^{ij61}+e^{ij34}) +a_{i6j}(e^{ij15}+e^{ij42}) \bigr)\\
&=&-2f^*\vol_M\wedge\bigl((a_{143}+a_{253})e^6+(a_{142}+a_{362})e^5 +(a_{251}+a_{361})e^4\bigr).
\end{eqnarray*} 
Now the three components in the last line are equal in nature. For instance,
\[ 2(a_{143}+a_{253})=\langle R^*_{e_1,e_3}U,e_4\rangle+ \langle R^*_{e_2,e_3}U,e_5\rangle=\Tr{R^\nabla_{\cdot,e_3}e_0}=r(e_3,e_0)\]
with $r$ the Ricci tensor of $M$. Therefore
\[ \dx*\phi=-f^*\vol_M\wedge\sum_{i=1}^3r(e_i,e_0)e^{i+3}=-f^*\vol_M\wedge\,(\ric U)^\flat \]
as we wished. Finally, the expression above also shows that $\dx*\phi=0$ if, and only if, $\ric U$ is a multiple of $U$. As it is well known (eg. from the decomposition of the curvature tensor), such multiple has to be a constant along the fibres $\inv{f}(x),\ x\in M$. Then this is also well known to imply $M$ is an Einstein manifold, ie. $\ric_x=s\,\Id,\ \forall x\in M$.
\end{proof}
We may also write $\delta\phi=-*\dx *\phi=-\ric U\lrcorner\alpha$.

The Riemannian manifold $S_M$ has a rather rich structure, since it is furnished with four \textit{volume} forms $\alpha,\alpha_1,\alpha_2,\mu\wedge\beta$ (not a differential system)  and one null-divergent unitary vector field $\theta^tU$. Notice the fibres of this bundle are all associative, ie. in each of them the restriction of $\phi$ is exactly the volume form. Notice furthermore that we may study other linear combinations with coefficients in $\cinf{M}$ of those four forms, rather than $\phi$, in order to obtain definite $G_2$ structures. Other problems relate with the exact 2-form $\beta$, if we ask when is it the K\"ahler form of a hermitian manifold transverse to the integral curves of the canonical vector field.

We remark also the interesting feature of $S_M$ which follows by changing the sign of the metric only on the fibre direction. This corresponds with the split octonions and the non-compact dual of $G_2$. Hence we may construct a $\widetilde{G}_2$ structure using the same method.

\subsection{The torsion forms for the constant sectional curvature case}
\label{ttfftcsc}

We may take from well known references the irreducible decomposition of the exterior algebra of $(\R^7)^*$ as a $G_2$ module (see eg. \cite{Bryant2,Fer}). The non trivial task resumes to degrees 2 and 3 since the star operator commutes with the group product. We have
\begin{equation}\label{decomp} 
\Lambda^2=\Lambda^2_{7}\oplus\Lambda^2_{14},\hspace{2cm} \Lambda^3=\Lambda^3_{1}\oplus\Lambda^3_{7}\oplus\Lambda^3_{27},
\end{equation}
where $\Lambda^2_{7}=\{\gamma\in\Lambda^2|\ \gamma\wedge\phi=-2*\gamma\}$, $\Lambda^2_{14}=\{\gamma\in\Lambda^2|\ \gamma\wedge\phi=*\gamma\}\simeq\g_2$, $\Lambda^3_1=\R\phi$, $\Lambda^3_7=\{*(\gamma\wedge\phi)|\ \gamma\in\Lambda^1\}$, $\Lambda^3_{27}=\{\gamma\in\Lambda^3|\ \gamma\wedge\phi=\gamma\wedge*\phi=0\}$. The indices below stand for the dimensions. Thence, there are isomorphic equivalents for degrees 4 and 5.

These irreducible summands induce corresponding subspaces $\Omega_j^p$, $j=1,7,14,27$ in the space of $p$-forms of a manifold with a $G_2$ structure and are used to classify such structures. The unique components $\tau_i$ of $\dx\phi$ and $\dx*\phi$ are called the torsion forms. One of them occurs in two places:
\begin{equation}\label{torsoes}
\dx\phi=\tau_0*\phi+3\tau_1\wedge\phi+*\tau_3,\hspace{1.7cm} \dx*\phi=4\tau_1\wedge*\phi+\tau_2\wedge\phi
\end{equation} 
with $\tau_i\in\Omega^i,\ \tau_2\in\Omega^2_{14},\ \tau_3\in\Omega^3_{27}$. Thus there are in principle sixteen classes of $G_2$ structures. 

In the case of $(S_M,\phi)$, we shall deduce the torsion forms in an article to follow, namely detecting scalar, Ricci and Weyl parts of $\cal R$ in the formula of $\dx\phi$. In particular we have found that $\tau_3$ never vanishes. For the moment one can read from theorem \ref{dphiemframe} more clearly, since those torsions have come up in a heavy manner. One evidence of this is the case of $M$ flat. Then $\dx*\phi=0$, due to proposition \ref{cocali}, and $\dx\phi=-\beta^2=\frac{6}{7}*\phi-*(*\beta^2+\frac{6}{7}\phi)$, showing which torsions do not vanish. We follow on just with the case of $M$ with constant sectional curvature.
\begin{prop}\label{tffcsc} 
If $M$ has constant sectional curvature $C$, then $\tau_1=\tau_2=0$ and
\begin{equation}
\tau_0=\frac{6}{7}(C+1),\hspace{1cm}\tau_3=(3C-\tau_0)\alpha+(2-\tau_0)\mu\wedge\beta-(C-\tau_0)\alpha_2.
\end{equation} 
\end{prop}
\begin{proof}
Since $M$ is Einstein, $\tau_1,\tau_2$ follow as in the flat case. Let the Riemann curvature tensor be $R^\nabla_{X,Y}Z=C(\langle Y,Z\rangle X-\langle X,Z\rangle Y)$ corresponding to constant sectional curvature $C$. Then it is easy to deduce that, in the notation of theorem \ref{dphiemframe}, ${\cal R}^{0i}=\langle R^\nabla\: e_0,e_i\rangle=-C e^{0i}$. Hence
\begin{eqnarray*}
\dx\phi &=& -C(e^{0156}+e^{0264}+e^{0345})-\beta^2+3C f^*\vol_M\\
&=& -C\mu\wedge\alpha_1-\beta^2+3C f^*\vol_M
\end{eqnarray*}
On the other hand $\phi\wedge\dx\phi=7\tau_0\vol_{S_M}$. Doing the same wedge with the result we found before, we deduce the equation 
\[ C\alpha_2\wedge\mu\wedge\alpha_1-\mu\wedge\beta^3+3C\vol_{S_M} =7\tau_0\vol_{S_M}.\]
This is equivalent to $3C+6+3C=7\tau_0$ and hence the values of $\tau_0$ and $\tau_3$ follow.
\end{proof}

\subsection{Some examples}
\label{se}

Now we apply the above to the most simple known examples of 4-manifolds.
The 4-torus with the flat metric gives a co-calibrated $S^3\times\T^4$ which seems to be unknown.

For $M=S^4$ we have $S_M=SO(5)/SO(3)$, since this is locally the product of two spheres. Otherwise we may start by proving that an isometry of the base space $M$ induces a $G_2$-isometry of $S_M$, ie. an isometry $g$ such that $g^*\phi=\phi$. This follows from the construction; as well as the conclusion that in the case of $S^4$ the action of $SO(5)$ is actually transitive. Then we have a co-calibrated 3-form such that
\begin{equation}\label{dphideS4} 
\dx\phi=*\phi+2\alpha+\mu\wedge\beta
\end{equation} 
differing from the case $\dx\psi=*4\psi$ of a well known nearly parallel $G_2$ structure $\psi$ on $SO(5)/SO(3)$, which may be seen in \cite{Bryant1,FriKaMoSe}. (A nearly parallel structure $\psi$ is one for which $\tau_0\in\R$ and all other torsions vanish. Such is the case too of the Hopf bundle $S^7\rr\Pro^1(\Hamil)$, which is an $SU(2)$-bundle given by the spin structure of the 4-sphere.)

For $M=\Pro^2(\C)=SU(3)/S(U(2)\times U(1))$ the action of $SU(3)$ also lifts to a tran\-sitive action on $S_M$ by $G_2$-isometries, since $T_zM=\C^3/z$; hence $S_M=SU(3)/U(1)$, an Aloff-Wallach space, where $w\in U(1)$ is included as the diagonal matrix $\mathrm{diag}(w,w,w^{-2})$. (Again this has a contrepart with the twistor space of $M$, the manifold of flags $SU(3)/\T^2$.) But the case of hermitian surfaces deserves to be studied in a proper place.

Recall no other compact irreducible symmetric spaces admit a transitive lift of the action to its sphere bundle (this follows from theorem 10.90 in \cite{Besse}, after M. Berger, and the list of those spaces which have rank 1).

It is also easy to deduce, from $\phi\wedge\dx\phi=7\tau_0\vol_{S_M}$, that
\begin{equation}\label{tau0} 
\tau_0=\frac{1}{7}(2r(U,U)+6)
\end{equation} 
in the general case. In particular, if $M$ is locally isometric to the standard hyperbolic 4-space, then $\tau_0=\tau_1=\tau_2=0$. Other developments in the general case have shown $\tau_3$ is never 0.

\end{document}